\documentclass[12pt,oneside]{amsart}
\usepackage[a4paper, total={6in, 9in}]{geometry}
\usepackage{cite,comment}
\usepackage{mathtools,amsmath}
\usepackage{comment}

\usepackage{amssymb}
\usepackage{layout}
\usepackage{tikz-cd}
\usepackage{amsmath,amsthm,mathrsfs}
\usepackage{marvosym}
\usepackage{xcolor}
\usepackage{colonequals}
\usepackage{enumitem}
\usepackage{mathdots}
\usepackage{framed}
\usepackage{url}
\usepackage{multicol}
\usepackage{fancyhdr}
\usepackage{nccmath}
\usepackage{color}
\usepackage{xparse}
\usepackage{bbm}
\usepackage[unicode,colorlinks]{hyperref}
\usepackage{microtype}
\usepackage{marginnote}
\usepackage{mathtools}
\usepackage{arydshln}
\usepackage{booktabs}

    \definecolor{dark-red}{rgb}{0.54,0,0}
    \definecolor{dark-green}{rgb}{0,0.54,0}
    \definecolor{dark-magenta}{rgb}{0.54,0,0.54}
    \definecolor{dark-cyan}{rgb}{0,0.54,0.54}
    \hypersetup{%
        unicode=true,
        breaklinks=true,
        pdfcreator=,
        pdfproducer=,
        pdfstartview =, 
        pdfpagelayout = ,
        pdfhighlight = /O,
        linkcolor=dark-red, 
        linkbordercolor=dark-red, 
        anchorcolor=red, 
        citecolor=dark-green, 
        citebordercolor=dark-green, 
        filecolor=dark-cyan, 
        filebordercolor=dark-cyan, 
        menucolor=dark-red, 
        menubordercolor=dark-red, 
        runcolor=dark-cyan, 
        runbordercolor=dark-cyan, 
        urlcolor=dark-cyan, 
        urlbordercolor=dark-cyan, 
}

\newcommand{\Af}{\mathbb{A}}

\newcommand\NN{\protect\mathbf{N}}
\newcommand\QQ{\protect\mathbf{Q}}
\newcommand\RR{\protect\mathbf{R}}
\newcommand\ZZ{\protect\mathbf{Z}}

\newcommand\PP{\protect\mathbb{P}}

\newcommand\cO{\mathcal{O}}

\DeclareMathOperator\ord{ord}

\theoremstyle{theorem} \newtheorem{proposition}{Proposition}[section]
\theoremstyle{definition} \newtheorem{definition}[proposition]{Definition}
\theoremstyle{theorem} \newtheorem{lemma}[proposition]{Lemma}
\theoremstyle{remark} \newtheorem{remark}[proposition]{Remark}
\theoremstyle{remark} 
\theoremstyle{theorem} \newtheorem{question}[proposition]{Question}
\theoremstyle{remark} 
\theoremstyle{definition} \newtheorem{example}[proposition]{Example}
\theoremstyle{definition} 
\theoremstyle{theorem} 
\theoremstyle{theorem} 
\theoremstyle{theorem} \newtheorem{theorem}[proposition]{Theorem}
\theoremstyle{theorem} 
\theoremstyle{definition} 
\theoremstyle{theorem} 
\theoremstyle{remark} 
\theoremstyle{definition} 
\theoremstyle{definition} 
\theoremstyle{definition} 
\theoremstyle{definition} 
\theoremstyle{remark} \newtheorem*{claim*}{Claim}
\theoremstyle{remark} 
\theoremstyle{theorem} 
\theoremstyle{theorem} 
\theoremstyle{definition} 
\theoremstyle{definition} 
\theoremstyle{theorem} 
\theoremstyle{remark}

\theoremstyle{definition} 
\theoremstyle{remark}
\theoremstyle{theorem} 
\theoremstyle{theorem} 

\DeclareMathOperator\Proj{Proj}

\DeclareMathOperator\lct{lct}
\DeclareMathOperator\Supp{Supp}
\DeclareMathOperator\mld{mld}
\DeclareMathOperator\codim{codim}

\DeclareMathOperator\CDiv{CDiv}
\DeclareMathOperator\Cosupp{Cosupp}
\DeclareMathOperator\Exc{Exc}
\newcommand\fa{\mathfrak{a}}
\newcommand\la{\lambda}
\newcommand\fb{\mathfrak{b}}
\newcommand\aeb{\mathfrak{a}^E_{\bullet}}
\newcommand\aem{\mathfrak{a}^E_{m}}
\newcommand\ab{\mathfrak{a}_{\bullet}}
\newcommand\fm{\mathfrak{m}}

\newcommand{\ld}[3]{a_{#1} \left( #2, #3 \right)}

{\endMakeFramed\noindent\ignorespacesafterend}

\numberwithin{equation}{section}

      \makeatletter
      \def\@setcopyright{}
      \def\serieslogo@{}
      \makeatother
   
\begin{document}

%



   \author{Harold Blum}
   \address{Department of Mathematics, University of Michigan, Ann Arbor MI 48105}
   \email{blum@umich.edu}



   

   \title[On Divisors Computing MLD's and LCT's]{On Divisors Computing MLD's and LCT's}

   \begin{abstract}
   
   We show that if a divisor centered over a point on a smooth surface computes a minimal log discrepancy, then the divisor also computes a log canonical threshold. 
 To prove the result, we study the asymptotic log canonical threshold of the graded sequence of ideals associated to a divisor over a variety. We systematically study this invariant
and also prove a result describing which divisors compute asymptotic log canonical thresholds. 
   \end{abstract}



   \thanks{This work was partially supported by NSF grant DMS-0943832.}




   \maketitle



   \section{Introduction}\label{intro}
   
   The \emph{log canonical threshold} and \emph{minimal log discrepancy} are two invariants of singularities that arise naturally in the study of birational geometry.
  Minimal log discrepancies are of particular interest due to a work of Shokoruv \cite{Sho} in which he proved that two conjectures on minimal log discrepancies (semicontinuity and the ascending chain condition (ACC)) imply the termination of flips, a result needed to complete the minimal model program in full generality.
  
  Shokurov originally conjectured that both the set of minimal log discrepancies and log canonical thresholds in fixed dimension should satisfy the ACC. The conjecture was proved for log canonical thresholds on smooth varieties  in \cite{FEM} and later in full generality \cite{MR3224718}. The general form of the ACC conjecture for minimal log discrepancies remains open. In this way, as well as others, minimal log discrepancies are less well understood than log canonical thresholds.

In order to to define these two invariants, we recall the following notion. A \emph{divisor over} a variety $X$ corresponds to a prime divisor on a normal variety $Y$, proper and birational over $X$.
 We call $(X, \fa^\la)$ a \emph{pair} if $X$ is a normal $\QQ$-Gorenstein variety, $\fa\subseteq \cO_X$ a nonzero ideal, and $\la\in \RR_{\geq 0}$. 
 The \emph{log discrepancy} of a pair $(X, \fa^\la)$ along $E$ is defined as 
 \[
 \ld{E}{X}{\fa^\la}  := k_E+1 - \la \ord_E(\fa)
, \]
where $k_E $ is the coefficient of $E$ in the relative canonical divisor and $\ord_E$ is the valuation given by order of vanishing along $E$.
A pair $(X,\fa^\la)$ is klt (resp., log canonical) if for all divisors $E$ over $X$, $a_E(X,\fa^\la)>0$ (resp., $\geq 0$). 
See Section \ref{ld} for further details on log discrepancies.

Arising from these definitions are two invariants that measure the ``nastiness'' of a singularity. Assuming $X$ has klt singularities, the \emph{log canonical threshold} of a nonzero ideal $\fa$ on $X$ is defined as 
\[
\lct(\fa) := \sup \{ \la \in \RR_{\geq0}\vert \, (X,\fa^\la) \text{ is log canonical} \}.\]
Given a klt pair $(X,\fa^\la)$ and a (not necessarily closed) point $\eta \in X$, the \emph{minimal log discrepancy} of $(X,\fa^\la)$ at $\eta$ is defined as
\[
\mld_\eta (X, \fa^\la) = \min \{ a_E(X,\fa^\la) \,\vert \,E \text{ a divisor over $X$ with $c_X(E) = \overline{\eta} $}  \}
.\]
See Section \ref{prelim} for further details on these definitions. 

In understanding these two invariants it is natural to make the following definition.
Given a divisor $E$ over $X$, we say that $E$ \emph{computes a log canonical threshold} if there exists a nonzero ideal $\fa$ on $X$ such that 
\[
\ld{E}{X}{\fa^ \la} = 0,\]
where $\la = \lct(\fa)$. Furthermore, we say that $E$ \emph{computes} $\lct(\fa)$.  
Similarly, we say that $E$ \emph{computes a minimal log discrepancy} if there exists a pair $(X,\fa^\la)$ such that 
\[
\mld_{\eta }(X,\fa^\la)
=
\ld{E}{X}{\fa^\la}
\]
with $\overline{\eta} = c_X(E)$.  Furthermore, we say that $E$ \emph{computes} $\mld_\eta (X, \fa^\la)$. 

\begin{question}
Which divisors over a variety compute log canonical thresholds (resp., minimal log discrepancies)?
\end{question}

Divisors computing log canonical thresholds satisfy special properties.
 As we will explain shortly, it is well known that divisors computing log canonical thresholds have finitely generated graded sequences of ideals. It is  not known if the same can be said for divisors computing minimal log discrepancies.

 It is trivial to see that  divisors computing log canonical thresholds also compute minimal log discrepancies. Indeed, if $E$ computes $\lambda = \lct(\fa)$ and $\overline{\eta} = c_X(E)$, then 
$E$ also computes $\mld_{\eta}(X, \fa^\la)$, which is $0$.
In the case of surfaces we prove the following.

 \begin{theorem}\label{surf}
 If $X$ is a smooth surface, then every divisor over $X$ centered at a point that computes a minimal log discrepancy also computes a log canonical threshold. 
 \end{theorem}

It has long been understood which divisors compute log canonical thresholds on smooth surfaces. 
For example, see \cite{KSC}, \cite{smith}, \cite{tucker}, and \cite{FJ}. In particular, 
\cite[Lemma 2.11]{FJ} implies that if
$E$ is a divisor over a smooth surface $X$ with $c_X(E)=\{x\}$ such that $E$ computes a log canonical threshold, then $\ord_E$ is a monomial valuation in some analytic coordinates at $x$.

Our methods extend to the case of du Val singularities. See Remark \ref{duvr} for our analogue of Theorem \ref{surf}.\\

Kawakita recently posted a related result which he proved independently\cite{Kawakita}. Kawakita
 showed that if $E$ is a divisor over a smooth surface $X$ with $c_X(E)=\{x\}$ such that $E$ computes a minimal log discrepancy, then $\ord_E$ is a monomial valuation in some local coordinates at $x$.\\

In proving Theorem \ref{surf}, we make use of the following object. Given a divisor $E$ over a normal variety $X$ and a morphism $f:Y \to X$ with $Y$ normal and $E$ arising as a prime divisor on $Y$, there is a corresponding graded sequence of ideals $\aeb=\{ \aem\}_{m\in \NN}$ defined as 
\[
\aem := f_\ast \cO_Y(-mE)
\]
(this is the ideal on $X$ of functions vanishing to order at least $m$ along $E$, hence it is independent of the choice of $f$).  Note that a graded sequence of ideals $\ab= \{ \fa_m \}_{m \in \NN}$ on $X$ is a sequence of ideals on $X$ such that $\fa_m \cdot \fa_n \subseteq \fa_{m+n}$ for all $m,n \in \NN$. 
We say that a graded sequence of ideals $\ab$ is finitely generated if  the graded ring $R(\ab):= \oplus_{m\in \NN} \ab$ is a finitely generated $\cO_X$-algebra. 

We now find a criterion to classify which divisors over $X$ compute log canonical thresholds.  As noted in \cite{JM}, if $E$ computes $\lct(\fb)$, then  
\[
\lct(\fb)
=
\lct \left(\fa^E_{\ord_E(\fb)} \right)
.\]
Furthermore, $E$ computes a log canonical threshold if and only if $E$ computes $\lct(\aem)$ for some $m\in \NN$. This criterion can be better phrased using $\lct(\ab^E)$, the \emph{asymptotic log canonical threshold} of the graded sequence of ideals $\ab^E$ (see Section \ref{lctsection}). We prove the following statement.

 \begin{theorem}\label{lct}
 Let $E$ be a divisor over a klt variety $X$. The following conditions are equivalent.
 \begin{enumerate}
 \item The divisor $E$ computes a log canonical threshold.
 \item The divisor $E$ computes an asymptotic log canonical threshold. 
 \item The equality $\lct(\aeb)=k_E+1$ holds.
 \end{enumerate}
 \end{theorem}

 The implications (1) $\Rightarrow$ (2) $\Rightarrow$ (3) appear in \cite{JM}. The implication (3) $\Rightarrow (1)$ requires a finite generation statement arising from the MMP. We find the equivalence between (1) and (2) to be somewhat surprising. 
  
  Rather than immediately proving Theorem \ref{surf}, we systematically study the invariant $\lct(\aeb)$. First, we note the following interpretation of the finite generation of $\aeb$. 
\begin{theorem}\label{tfg}(\cite[Corollary 3.3]{Ishii})
Let $X$ be a normal variety and $E$ a divisor over $X$ such that $\codim_X (c_X(E)) \geq 2$. If $\aeb$ is finitely generated, then 
\[
Y:=\Proj_X( \oplus_{m\geq 0} \fa_m^E)  \to X\]
is a proper birational morphism of normal varieties with exactly one exceptional divisor $E_Y$. Furthermore, $E_Y$ and $E$ induce the same valuation on $K(X)$. ß
\end{theorem}

Such birational morphisms with exactly one exceptional divisor were studied in \cite{Ishii} and referred to as \emph{prime blowups}. In \cite{plt}, the author looked at \emph{plt-blowups} which can be thought of as prime blowups with mild singularities (see Definition \ref{pltdef}).  
The following proposition relates the value of $\lct(\aeb)$ to the model $\Proj_X( \oplus_{m\geq 0} \fa_m^E)  $.

\begin{proposition} \label{proplctgeom}
If $X$ is a klt variety and $E$ a divisor over $X$  with \mbox{$k_E < \lct(\aeb)$}, then the following hold.
\begin{enumerate}
\item The graded sequence $\aeb$ is finitely generated.
\item The map $Y:= \Proj_X( \oplus_{m\geq 0} \fa_m^E) \to X$ is a proper birational morphism of normal varieties. Additionally, there is a prime  divisor $E_Y \subseteq Y$ such that $-E_Y$ is relatively ample $X$ and $E_Y$ induces the same valuation as $E$ on $K(X)$.
\item The variety $Y$ is klt, and $\lct(Y,E_Y) = \lct(\aeb)-k_E$.
\end{enumerate}
\end{proposition}

 The first assertion of the above proposition
is an elegant restatement of the well known fact that if a divisor $E$  has log discrepancy in the interval [0,1) along a klt pair, then $\aeb$ is finitely generated. (When $X$ is $\QQ$-factorial, this statement follows from \cite[Corollary 1.39]{K13}.)

In the literature, there are two examples  of divisors over smooth varieties with non-finitely generated graded sequences of ideals \cite{CGP}  \cite{Kur}. Again, by \cite[Corollary 1.39]{K13}, it is well known that these divisors cannot compute log canonical thresholds (but not whether they compute minimal log discrepancies). As we explain in Example \ref{nfgex}, our results in Section \ref{sectionmld} show that the divisor in  \cite{Kur} does not compute a minimal log discrepancy.
\\

\noindent{\bf Structure of the Paper:} In Section \ref{prelim} we provide preliminary information on log discrepancies, log canonical thresholds, and graded sequence of ideals. Section \ref{sectionrees}  provides a proof of Theorem \ref{tfg} and a related statement of independent interest on the graded sequence of ideals corresponding to a divisor over a variety.
In Section \ref{sectionmmp} we prove Proposition \ref{proplctgeom}. Section \ref{sectionmld} gives of proof Proposition \ref{lct} and a condition on which divisors compute minimal log discrepancies. 
Section \ref{sectionsurf} concerns minimal log discrepancies on surfaces and gives a proof of Theorem \ref{surf}. Lastly, Section  \ref{examples} provides computations of $\lct(\aeb)$ in a few examples. 
  \\
  
\noindent{\bf Acknowledgements:}
I would like to thank my advisor Mircea Musta{\c{t}}{\u{a}} for introducing me to many of the topics discussed in this paper and guiding my research. I would also like to thank Mattias Jonsson and Karen Smith for useful conversations. 
 
 \section{Preliminaries}\label{prelim}
 
 \noindent {\bf Conventions:} For the purpose of this paper, a variety is a reduced, irreducible, separated
scheme of finite type over a field $k$. Furthermore, we will assume that $k$ is of characteristic 0 and algebraically closed.
 
 \subsection{Log Resolutions and Divisors Over a Variety}\label{divover}
 
 Let $X$ be a variety and $\fa\subseteq \cO_X$ a nonzero ideal. A morphism $f:Y\to X$ is a \emph{log resolution}
 of $(X, \fa)$ if $f$ is a projective birational morphism, $Y$ is smooth, $\Exc(f)$ has pure codimension 1, $\fa\cdot \cO_Y=\cO_Y(-D)$ for some effective divisor $D$ on $Y$, and $\text{Exc}(f) + D$ is a simple normal crossing divisor.
 
Let $f:Y\to X$ be a proper birational morphism of normal varieties, and $E$ a prime divisor on $Y$.  We define the \emph{center} of $E$ on $X$ to be $c_X(E)= f(E)$.  Arising from $E$ is a discrete valuation of $K(X)$ that we denote by $\ord_E$. The valuation corresponds to the DVR $\cO_{Y,E} \subseteq K(Y) \simeq K(X)$. Given a nonzero ideal $\fa$ on $X$, we set $\ord_E(\fa) = e$ where $\fa \cdot \cO_{Y,E}=(t^e)$ and $t\in \cO_{Y,E}$ is a uniformizing parameter of the DVR. We will use the convention that $\ord_E( {\bf 0} ) = +\infty$, where ${\bf 0}$ denotes the zero ideal.  For a $\QQ$-Cartier divisor $D\geq 0$, the vanishing of $D$ along $E$ is defined in the following way. Choose $m\in \ZZ_{>0}$ such that $mD$ is Cartier and let 
 \[
 \ord_E(D) \coloneqq  \frac{1}{m} \ord_E ( \cO_X(-mD))
 \,.\]
 By linearity, $\ord_E$ extends to a function on $\text{CDiv}_\RR(X)$. 
 
 Given $Y \to X$ and $Y' \to X$ as above, we identify a prime divisor on $Y$ with a prime divisor on $Y'$ if the two divisors induce the same valuation on $K(X)$. Note that if $Y'\to X$ factors as $Y' \to Y \to X$ and $E\subset Y$ and $E'\subset Y'$ prime divisors, then we identify $E$ and $E'$ if and only if $E'$ is the strict transform of $E$ on $Y'$. If $E$ and $E'$ are identified then $c_X(E)=c_X(E')$. A \emph{divisor over} $X$ is an equivalence class given by this relation. 
 
 \begin{remark}\label{remideal}
 Let $f:Y \to X$ be a proper birational morphism of normal varieties. We will often consider ideals on $X$ of the form
 \[
 \fa = f_\ast \cO_X(-a_1E_1 - \cdots - a_r E_r) \subseteq \cO_X \]
 where $E_i$ is a prime divisor on $Y$ and $a_i>0$ for $1 \leq i \leq r$. 
For an open affine $U \subseteq X$,  the ideal $\fa(U)$ can be understood as followed. 
If $U$ does not intersect $f(E_i)$ for all $1\leq i \leq r$, then $\fa(U)= \cO_X(U)$. 
Otherwise, 
\[
\fa (U)  = \left\{ f\in \cO_X(U) \vert \ord_{E_i}(f) \geq a_i \text{ for } i \text{ such that } c_X(E_i)\cap U \neq \emptyset \right\}.\]
  \end{remark}
 
 \subsection{Rees Valuations}
For a nonzero ideal $\fa$ on a normal variety $X$,  let $Z=\widetilde{B_{\fa}(X)} \overset{\pi}{\to} X$ denote the normalized blowup of $X$ along $\fa$.
 Note that $\fa\cdot \cO_Z = \cO_{Z}(-D)$ for some effective divisor $D$ on $Z$. 
 The \emph{integral closure of $\fa$} is the ideal
 \[
 \overline{\fa} := \pi_\ast \cO_Z (-D)
 \,.\]
It is always the case that $\fa \subseteq \overline{ \fa}$.  
 We say that $\fa$ is \emph{integrally closed} if $\fa = \overline{\fa}$.
The \emph{Rees valuations of $\fa$} are the valuations of $K(X)$ corresponding to prime divisors in the support of $D$. The following propositions provide information on Rees valuations and integral closures of ideals.

\begin{proposition}\label{proprees}
Let $\fa$ be a nonzero ideal on a normal variety $X$. The set of Rees valuations of $\fa$ are the valuations corresponding to the smallest set of divisors $\{E_1,\ldots, E_r\}$ over $X$ such that
\[
\overline{\fa^m} = \{ f\in  \cO_X \vert \ord_{E_i}(f) \geq  m\cdot \ord_{E_i} (\fa) \text{ for } 1\leq i \leq r \}
\]
for all $m\in \ZZ_{>0}$.
\end{proposition}
\begin{proof}
See\cite[Section 10.2]{Swanson}.
\end{proof}

 \begin{proposition}\label{reesdet}
  Let $\fa$ be a nonzero ideal on a normal variety $X$. 
The ideal $\fa$ is integrally closed if and only if there exists a proper birational morphism $f:Y\to X$ with $Y$ normal and $D$ an effective cartier divisor on $Y$ such that
 \[
 \fa = f_\ast \cO_Y(-D) 
 .\]
 \end{proposition}
 
 \begin{proof}
The forward direction follows immediately from the definition of a Rees valuation.
See \cite[Proposition 9.6.11]{MR2095472} for the reverse direction.
 \end{proof}
 
 \begin{proposition}\label{blowupnormal}
If  $X$ is a normal variety and $\fa$ a nonzero ideal on $X$ such that $\fa^n$ is integrally closed for all $n\in \ZZ_{>0}$, then the blowup of $X$ along $\fa$ is normal. 
 \end{proposition}
 \begin{proof}
 This follows from \cite[Proposition 5.2.1]{Swanson}
  \end{proof}
  \subsection{Log Discrepancies} \label{ld}
 
 Let $X$ be a normal variety. We say that $X$ is $\QQ$-Gorenstein if the canonical divisor $K_X$ is $\QQ$-Cartier. 
If $f:Y\to X$ is a proper birational morphism of normal varieties and $X$ is $\QQ$-Gorenstein, the relative canonical divisor of the morphism is defined as $K_{Y/X}:= K_Y - f^\ast K_X$ where $K_Y,K_X$ are chosen so that $f_\ast K_Y = K_X$. 

   If $E$ is a divisor over $X$ that appears as a prime divisor on $Y$, we let $k_E$ denote the coefficient of $E$ in $K_{Y/X}$.  When the base variety is unclear, we will use the notation $k_{E,X}$. The value $k_E$ is not dependent on the model $Y$ but on the valuation $\ord_E$.  Thus, this invariant extends to an invariant for a divisor over $X$.

 Given a  normal, $\QQ$-Gorenstein variety $X$, a nonzero ideal $\fa\subseteq \cO_X$, and $\la \in \RR_{\geq 0}$, we refer to $(X,\fa^\la)$ as a \emph{pair} and define the \emph{log discrepancy} of $(X, \fa^\la)$ along $E$ as
 \[
 \ld{E}{X}{ \fa^\la} \coloneqq  k_E+1-\la \ord_E(\fa).\]
 We say that the pair $(X, \fa^\la)$ has klt (resp. log canonical) singularities if $\ld{E}{X}{\fa^\la}>0$ (resp. $\geq0)$ for all divisors $E$ over $X$.
 
 If $X$ is $\QQ$-Gorenstein and $\Delta$ is $\QQ$-factorial, the log discrepancy of $(X,\Delta)$ along $E$ is 
 \[
 \ld{E}{X}{ \Delta} : = k_E+1- \ord_E( \Delta).\]
 We say that $(X,\Delta)$ is a klt (resp. log canonical) pair if  $\ld{E}{X}{ \Delta}>0$ (resp. $\geq0$) for all divisors $E$ over $X$.  
  We say that $(X,\Delta)$ is plt if $\ld{E}{X}{ \Delta}>0$ for all divisors $E$ over $X$ with $\codim c_X(E) \geq 2$. 
 Lastly, we say that $X$ is klt (resp. log canonical) if $X$ is normal, $\QQ$-Gorenstein, and $(X,0)$ is klt (resp. log canonical). 
 
  \subsection{Minimal Log Discrepancies.} Given a pair $(X,\fa^\la)$ and a (not necessarily closed) point $\eta \in X$, we define 
 the \emph{minimal log discrepancy} of $(X,\fa^\la)$ at $\eta$ to be 
 \[
 \mld_\eta(X,\fa^\la)\coloneqq \inf  \left\{   \ld{E}{X}{\fa^\la} \vert \text{$E$ is a divisor over $X$ with $c_X(E) = \overline{ \{ \eta \}}$} \right\}
 .\]
Assuming $\codim(\overline{ \{\eta \} }) \geq 2$,  the value of $\mld_\eta(\fa^\la)$ is either $\geq0$ or $=-\infty$. If the above infimumum is $\geq 0$, it is necessarily a minimum. See \cite{Ambro} for more details on minimal log discrepancies. 
 \subsection{Graded Sequences of Ideals}\label{gs}
 
 A \emph{graded sequence of ideals} on a variety $X$ is a sequence of ideals $\ab = \{ \fa_m\}_{m\in \NN }$ on $X$ 
 such that $\fa_m \cdot \fa_n \subseteq \fa_{m+n}$ for all $m,n\in\NN$. (We use the convention that $\NN$ includes $0$.) In the case when $\fa_0 = \cO_X$, the \emph{Rees algebra} of
 $\ab$ is the $\NN$-graded $\cO_X$-algebra
 \[
 R(\ab) = \bigoplus_{m\in \NN}\fa_m
 .\]
 We say that $\ab$ is \emph{finitely generated} if $R(\ab)$ is a finitely generated $\cO_X$-algebra. 
 We say that $\ab$ is \emph{finitely generated in degree $n$} if $\fa_{nm}= (\fa_n)^m$ for all $m\in \NN$. A graded sequence  of ideals $\ab$ is finitely generated if and only if $\ab$ is finitely generated in some degree $n\in \ZZ_{>0}$ \cite[Lemma 2.1.6.v]{ega}.
 
 We list three examples of graded sequences of ideals that arise in algebraic geometry. 
 \begin{itemize}
 \item For a trivial example, let $\mathfrak{b}$ be an ideal on $X$ and define $\mathfrak{a}_m= \mathfrak{b}^m$. 
 \item 
 
Let $f:Y\to X$ be a proper birational morphism of normal varieties and $E$ be a prime divisor on $Y$. The divisor $E$ gives rise to a graded sequence of ideals, denoted by $\aeb$, defined by
\[
\fa_m^E \coloneqq  f_\ast \cO_Y(-mE)
.\] Note that this only depends on $E$
and not on the model $Y$.
\item Let $\mathcal{L}$ be a line bundle on a variety $X$, having nonnegative Kodaira dimension and $\fa_m(\mathcal{L})$ denote the base locus of $\left| \mathcal{L}^m \right|$. This example was studied in \cite{base}.
\end{itemize}

For a graded sequence $\ab$, let $S(\ab) = \{ m\in \NN \vert \fa_m \neq {\bf 0)} \}$, where ${\bf 0}$ is the zero ideal. 
If $\ab$ is a graded sequence of ideals on $X$ such that $S(\ab)$ is nonempty, we define
\[ \ord_E(\ab)\coloneqq  \lim_{m\to \infty , m \in S(\ab) } \frac{1}{m} \ord_E (\fa_m)= \inf_{m\geq 1}  \frac{1}{m} \ord_E (\fa_m).\]
See \cite[Section 2]{JM} for further details. 

\subsection{Log Canonical Thresholds.} \label{lctsection} For a nonzero ideal $\fa$ on a klt variety $X$, the \emph{log canonical threshold} of the ideal is defined as
\[
\lct(\fa )\coloneqq  \sup\{ \lambda \in \RR_{\geq 0} \vert (X, \fa^\lambda) \text{ is log canonical} \}.\] From this definition, it follows that $\lct( \cO_X)= +\infty$. We define $\lct( {\bf 0} ) = 0$, where ${\bf 0}$ denotes the zero ideal. Note that $(X, \fa^\la)$ is log canonical if and only if 
\[
k_E+1 - \la \ord_E(\fa ) \geq 0\]
for all $E$ over $X$. 
Thus, 
\[
\lct( \fa) =  \inf \left\{ \frac{k_E+1}{ \ord_E(\fa)} \vert \text{\,$E$ a divisor over $X$ with $\ord_E( \fa)>0$  } \right\}.\]
Moreover, if $Y\to X$ is a log resolution of $(X, \fa)$ and $\fa \cdot \cO_Y= \cO_Y(-D)$,
it is sufficient to take the above infimum over all prime divisors contained in $\Supp(D)$ (see \cite[Section 9.3]{MR2095472}). Since $\Supp(D)$ has finitely many components, the above infimum is necessarily a minimum. 
Similarly, if $\Delta \in \CDiv_{\RR\geq 0} (X)$, then the log canonical threhsold of $\Delta$ is
\[
\lct(X,\Delta) := \sup \{ \la \in \RR_{\geq 0} \vert (X,\lambda \Delta) \text{ is log canonical} \} .\]

In Section \ref{intro}, we said that for a nonzero ideal $\fa$, a divisor $E$ computes $\lct(\fa)$ if 
$\ld{E}{X}{ \fa^{ \la}}=0$,
where $\la = \lct(\fa)$. 
Note that this is equivalent to the condition that $\lct(\fa)= (k_E+1)/ \ord_E( \fa)$. 
 
The following elementary properties of log canonical thresholds give rise to the definition of the log canonical threshold for a graded sequence of ideals. 

\begin{lemma}\label{easy}
Let $\fa,\fb$ be nonzero ideals on a klt variety $X$. The following hold:
\begin{enumerate}

\item $\lct(\fa)= m \cdot \lct(\fa^m)$.

\item $\lct(\fa) \leq \lct(\fb)$ if $\fa \subseteq \fb$.

\end{enumerate}
\end{lemma}
Thus, if $\ab$ is a graded sequence of ideals, we have that $m\cdot \lct(\fa_m) \leq mn \cdot \lct(\fa_{mn})$. If $S(\ab)$ is nonempty, we define the \emph{asymptotic log canonical threshold} of $\ab$ to be
\[
\lct(\ab)\coloneqq  \lim_{m\to \infty, m \in S(\ab)} m\cdot \lct(\fa_m)= \sup_{m\geq 1} m\cdot \lct(\fa_m).\]
For the second equality, see \cite[Lemma 2.3]{JM}. 

 In the following statements we collect some basic information on this asymptotic invariant. The content of the following statement is proven in \cite{JM}.

 \begin{proposition} \label{tlct}
Let $X$ be a klt variety and $\ab$ a graded sequence of ideals on $X$ such that $S(\ab)$ is nonempty. The following hold:
\begin{enumerate}
\item $\lct(\ab)=\inf_{F} \frac{k_F+1}{\ord_F(\ab)}$, where the infimum runs over all divisors $F$
over $X$. 

\item If $\ab$ is finitely generated in degree $m$ (i.e. $\fa_{mn} = (\fa_m)^n$ for all $n\in \NN$), then 
$\lct(\ab)= m\cdot \lct(\fa_m)$. Furthermore, if $f:Y\to X$ is a log resolution of $(X, \fa_m)$ with $\fa_m \cdot \cO_Y=\cO_Y(-D)$, then
\[
\lct(\ab) = \min_{E \subseteq \Supp(D)}  \left\{  \frac{k_E+1}{\ord_E(\ab)} \right\}
.\]

\end{enumerate} 

\end{proposition}

We say that $E$ computes the \emph{asymptotic log canonical threshold} of $\ab$ if $\lct(\ab)=(k_E+1)/\ord_E(\ab)$

\begin{proof}
See \cite[Corollary 2.16]{JM} for (1). The proof of (2) follows from the definition of the log canonical threshold of a graded sequence of ideals.
\end{proof}

 \begin{lemma}\label{self}
Let $E$ be a divisor over a variety $X$. We have that 
\begin{enumerate}
\item $\ord_E(\fa_m^E) = m$ for all $m \in \ZZ_{>0}$ divisible enough,
\item $\ord_E(\fa_m^E)=m$ if $\aeb$ is finitely generated in degree $m$, and
\item $\ord_E(\aeb)=1$, which implies that $\lct(\aeb) \leq k_E+1$.  
\end{enumerate}
 \end{lemma}
 \begin{proof}
 Since elements of $\fa_m^E$ vanish to at least order $m$ along $E$, we see that $\ord_E(\fa_m^E) \geq m$. Set $n= \ord_E(\fa_1^E)$. Since $(\fa_1^E)^m \subset \fa_{n \cdot m}^E$, 
 we see that $\ord_E(\fa_{n\cdot m} ) \leq n\cdot m$. Thus, (1) follows.  

Next, (2) follows from (1) after noting that if $\aeb$ is finitely generated in degree $m$, then 
\[
\ord_E(\fa_{m\cdot n} ^E) = \ord_E( (\fa_m^E) ^n) = n \ord_E(\fa_m^E).\] Finally, (3) follows immediately from (1) and Proposition \ref{tlct}. 
 \end{proof}

The above lemma gives that $\lct(\aeb) \leq k_E+1$. The following proposition explains the difference between these values. 

\begin{proposition} \label{lctexp}
For $X$ a klt variety and $E$ a divisor over $X$, we have that 
\[
 k_E+1- \lct(\aeb)
 =
\inf \left\{ a_E(X, \fa^\la) \vert \, (X,\fa^\la) \text{ is a log canonical pair} \right\}  .\]
\end{proposition}

\begin{proof}
We first note that 
\[
k_E+1 - \lct(\aeb) = \inf_{m\geq 1}  \{ k_E+1 - m \cdot \lct(\fa_m^E)  \} =  \lim_{m\to \infty}   \left( k_E+1 - m\cdot \lct(\fa_m^E)  \right)
. \]
By Lemma \ref{self}, for $m$ divisible enough $\ord_E(\fa_m^E)= m$, and, thus,
\[
k_E+1 -m \cdot \lct(\fa_m^E) = \ld{E}{X}{(\fa_m^E)^{\lct(\fa_m^E)} } 
.\]
Since $(X,(\fa_m^E)^{\lct(\fa_m^E)})$ is a log canonical pair,
the relation $\geq$ of the desired equality follows. 

Next, we show the $\leq$ relation of the equality.
Let $(X, \fb^\la)$ be a log canonical pair and set $m' = \ord_E(\fb)$ and $\la' = \lct(\fa_{m'}^E)$. 
Since $\fb \subseteq \fa_{m'}^E$, we have that $(X, \fa_{m'}^\la)$
is log canonical as well, and it follows that $\la \leq  \la'$. Thus, 
\[
\ld{E}{X}{ \fb^\la}   \geq \ld{E}{X}{ (\fa_{m'}^E)^\la} \geq \ld{E}{X}{ (\fa_{m'}^E)^{\la'} }
\geq 
k_E+1 - m' \cdot \lct(\fa_{m'}^E)
,
\]
where the first inequality follows from the relation $\fb \subseteq \fa_{m'}^E$, the second from the relation $\la \leq \la'$, and the last inequality from $m' \leq \ord_E( \fa_{m'}^E)$. The result follows. 
\end{proof}
\section{Rees Valuations and Graded Sequences}\label{sectionrees}
In this section, we prove Theorem \ref{tfg}. While the statement was proven in \cite 
{Ishii}, we give a short proof since it will be useful for us to understand the proof. We begin with the following lemma. 
\begin{lemma}\label{lemrees} Let $E$ be a divisor over a normal variety $X$ such that $\aeb$ is finitely generated. If $\aeb$ is finitely generated in degree $m$, then the following hold:
\begin{enumerate}
\item The positive powers of the ideal $\fa_m^E$ are integrally closed. 
\item The ideal $\fa_m^E$ has exactly one Rees valuation, namely $E$. 
\end{enumerate}
\end{lemma}

\begin{proof}
Let $f:Y\to X$ be a projective birational morphism such that $Y$ is normal, $E$ is a prime divisor on $Y$, and $\aeb$ is finitely generated in degree $n$. Then 
\[
f_\ast \cO_Y(-(n\cdot m)E) = \fa_{m \cdot n}^E = (\fa_n^E)^m,
\]
where the first equality comes from the definition of $\aeb$ and the second from our assumption that
$\aeb$ is finitely generated in degree $n$.

By Proposition \ref{reesdet}, the above equality gives that $(\fa_n^E)^m$ is integrally closed. By Proposition \ref{proprees}, the valuation corresponding to $E$ is the unique Rees valuation of $\fa_n^E$. \end{proof}

\begin{proposition}\label{propblowup}
If $E$ is a divisor over a normal variety $X$ such that $\aeb$ is finitely generated, then the following hold:
\begin{enumerate}
\item The map $Y:=\Proj_X( \oplus_{m\geq 0} \fa_m^E) \to X$ is a proper birational morphism of normal varieties. 
\item There exists a prime divisor $E_Y \subset Y$ that induces the same valuation as $E$ on $K(X)$. Furthermore, $E_Y$ is $\QQ$-Cartier and $-E_Y$ is relatively ample over $X$. 
\end{enumerate}
\end{proposition}
\begin{proof}
Since $\aeb$ is assumed to be finitely generated, we may choose $n\in \NN$ so that $\aeb$ is finitely generated in degree $n$. 
We claim that $Y$ is isomorphic to the blowup of $X$ along $\fa_n^E$. Indeed, $Y \simeq \Proj(\oplus_{m\geq 0} \fa_{nm}^E)$. Since $\aeb$ is finitely generated in degree $n$, $\Proj(\oplus_{m\geq 0} \fa_{nm}^E)= \Proj(\oplus_{m\geq 0} (\fa_{n}^E)^m)$. Thus, the claim is complete. 

Now, we apply the previous lemma. 
By Lemma \ref{lemrees} (1) combined with Proposition \ref{blowupnormal},  the variety $Y$ is normal. Since $\ord_E$ is unique Rees valuation of $\fa_n^E$ (Lemma 
\ref{lemrees} (2)), we see $\fa_n^E \cdot \cO_Y = \cO_Y(-nE_Y)$, where $E_Y$ is a prime divisor on $Y$ and $\ord_{E_Y} = \ord_E$. By the equality  $\fa_n^E \cdot \cO_Y = \cO_Y(-nE_Y)$, the divisor $E_Y$ is $\QQ$-Cartier and $-E_Y$ is relatively ample over $X$. This completes the proof.  
\end{proof}

\begin{proof}[Proof of Theorem \ref{tfg}]
The theorem is a consequence of the previous proposition. Note that $E_Y \subset Y$ is exceptional over $X$, since $\codim_X (c_X(E_Y))\geq 2$.ß
\end{proof}

If $E$ is a divisor over a normal variety $X$, it is natural to ask which exceptional divisors arise on the normalized blow up of $\fa_m^E$ for very divisible $m\in\NN$. Lemma \ref{lemrees}  says that if $\aeb$ is finitely generated, there is one such exceptional divisor, namely $E$. The following proposition says that in the general case (when $\aeb$ is not necessarily finitely generated) $E$ does arise on such a model.The proof of the proposition is similar in spirit to a proof of Artin (see \cite[Lemma 2.45]{KM}).

\begin{proposition}
Let $E$ be a divisor over a variety $X$. For $m$ divisible enough, $\ord_E$ is a Rees valuation of $\fa_m^E$. 
\end{proposition}

\begin{proof}
We assume that $X$ is affine (since blowups can be computed locally), and let $(R, \fm_R)\subset K(X)$ denote the DVR of the function field of $X$ associated to the divisor $E$. Thus, 
\[
\fa_m^E = \fm_R^m \cap \cO_X(X).\]
We fix the following notation: let $X_m$ be the blowup of $X$ along $\fa_m^E$, $x_m$ the image of the 
generic point of $E$ on $X_m$, and $\cO_m=\cO_{X_m,x_m}$. 

We must show that $\overline{ \{x_m \}}$ is codimension 1 on $\widetilde{X_m}$ (the normalization of $X_m$) for all $m$ divisible enough. In order to do this, it is sufficient to show that $\cO_m=R$
for $m\in \ZZ_{>0}$ divisible enough. \\

\noindent {\bf Claim 1}: For $m,n\in \ZZ_{>0}$ such that $\ord_E (\fa^E_m)=m$, we have $\cO_m \subseteq \cO_{mn}$.

Let $\fa_m^E = (z_1,\ldots,z_r)$ with $\ord_E(z_1)=m$. We have $z_1^n \in \fa_{mn}^E$. Thus, we
may choose elements $y_i \in \cO_X(X)$ so that $\fa_{mn} = (z_1^n, y_1,\ldots,y_s)$. The two rings
\[
\cO'_{m} \coloneqq  \cO_X(X) \left[ \frac{z_2}{z_1},\ldots, \frac{z_r}{z_1} \right] \hspace{0.4 in} \text{ and } \hspace{0.4 in} 
\cO'_{mn} = \cO_X(X) \left[ \frac{y_1}{z_1^n},\ldots, \frac{y_s}{z_1^n} \right]
\]
correspond to charts of the blowup of $X$ along $\fa_m^E$ and $\fa_{mn}^E$, respectively. 
Since the map $\cO_X(X) \to R$ extends to maps from the above rings (we are using that $\ord_E(z_i/z_1), \ord_E(y_i/z_1^n) \geq 0$), we see that $x_m$ and $x_{mn}$ lie on these charts. 
Thus, $\cO_m$ and $\cO_{mn}$ correspond to $\cO'_m$ and $\cO'_{mn}$ localized at their intersections with $\fm_R$. 

To show the desired inclusion, it is sufficient to show that $\cO'_m\subseteq \cO'_{mn}$ 
(i.e. $z_i/z_1$ lies in $\cO'_{mn}$ for all $1 \leq i \leq r$). We rewrite the fraction in the following way
\[
\frac{z_i}{z_1} = \frac{z_i z_1^{n-1}}{z_1^n} .\]
Since $\ord_E(z_i z_1^{n-1}) \geq \ord_E(z_1^n) \geq mn$, it follows that $z_iz_1^{n-1} \in \fa_{mn}^E$. Thus, 
$z_iz_1^{n-1}/z_1^n \in \cO'_{mn}$. 
\\

\noindent {\bf Claim 2}: For any element $f\in R$, there exists $m_f \in \ZZ_{>0}$ so that $f\in \cO_{m_f}$. 

Write $f=u/v$ where $u,v\in \cO_X(X)$. Set $m_f\coloneqq  \ord_E(v)$. Since $f\in R$, we have that
\[
m_f = \ord_E(v) \leq \ord_E(u)
,\]
and we can write $\fa_{m_f} = (u,v,w_1,\ldots,w_t)$ for some choice of $w_i$'s. Since $\ord_E(w_i/v)\geq 0$, $x_{m_f}$ lies on the chart of the blowup of $\fa_{m_f}$ corresponding to 
\[
\cO_0 \left[ \frac{u}{v} , \frac{w_1}{v},\ldots, \frac{w_t}{v} \right]
\]
and $\cO_{m_f}$ corresponds to the above ring localized at its intersection with $\fm_R$. Thus, $f\in \cO_{m_f}$. \\

The above two claims combine to show that there is an ascending sequence of subrings of $R$
\[
\cO_{1!} \subseteq \cO_{2!} \subseteq \cO_{3!} \subseteq \cdots
\]
whose union is $R$. As argued in the proof of \cite[Lemma 2.45]{KM}, for $m$ divisible enough $\cO_m=R$. 
This completes the proof. 
\end{proof}

\section{Finite Generation Using MMP}\label{sectionmmp}

Utilizing the relative LMMP (log minimal model program), we obtain a sufficient condition for the finite 
generation of $\aeb$, the graded sequence corresponding to a divisor $E$ over a variety $X$. We first state
the following proposition, which is known to experts. In the case when $X$ is $\QQ$-factorial, the proposition follows directly from  \cite[Corollary 1.39]{K13}. Our proof relies on the finite generation statement in \cite{BCHM}. 

\begin{proposition}\label{mldk}
Let $E$ be a divisor over a klt variety $X$. If there 
exists a divisor $\Delta$ such that $(X,\Delta)$ is klt and $a(E,X,\Delta)<1$, then $\aeb$ is finitely generated.
\end{proposition} 

Before proving the proposition, we recall the following. We say that two $\QQ$-divisors $D$ and $D'$ are \emph{$\QQ$-linearly equivalent}, denoted  by $D \sim_\QQ D'$, if there exists $m\in \ZZ_{>0}$ so that $mD$ and $mD'$ are linearly equivalent integral divisors. Fix a proper morphism $g:X\to Z$. We say that $D$ and $D'$ are \emph{$g$-linearly equivalent}, denoted by $D\sim_{\QQ,g} D'$ if $D-D'$ is a $\QQ$-linear combination of principal divisors and Cartier divisors pulled back from $Z$.

\begin{proof}
Let $(X,\Delta)$ be a klt pair and set $a:= \ld{E}{X}{\Delta}<1$. We may choose $f:Y\to X$ to be a log resolution of $(X, \Delta)$ such that $E_Y$ is a divisor on $Y$ identified with $E$. We have that 
\[
K_Y+f_\ast^{-1} \Delta  \sim_{\QQ} f^\ast(K_X +\Delta)+  (a-1)E_Y + \sum_{i=1}^{r} (\ld{E_i}{X}{\Delta} -1) E_i ,
\]
where $f_\ast^{-1}\Delta$ denotes the strict transform of $\Delta$ and $E_Y, E_1,\ldots, E_r$ are the exceptional divisors of $f$. Adding $(1-\epsilon) \sum E_i$ to both sides and looking at $f$-linear equivalence, we have 
\[
K_Y+f_\ast^{-1} \Delta  + (1-\epsilon) \sum E_i \sim_{\QQ ,f} (a-1) E_Y+
  \sum_{i=1}^{r} (\ld{E_i}{X}{\Delta} -\epsilon ) E_i  .\]
    
  Since $(X,\Delta)$ is klt, $\ld{E_i}{X}{\Delta} >0$ for all $i$. Thus, we may fix $0< \epsilon <1$ so that  $\ld{E_i}{X}{\Delta} -\epsilon >0$ for all $i$. Now, $(Y, f_\ast^{-1} \Delta  + (1-\epsilon) \sum E_i)$ is a klt pair, since $Y$ is smooth and  $f_\ast^{-1} \Delta  + (1-\epsilon) \sum E_i$ is a simple normal crossing divisor with coefficients in $[0,1)$. Additionally, $K_Y+  f_\ast^{-1} \Delta  + (1-\epsilon) \sum E_i$ is $f$-big. (Since $f$ is a birational morphism, all $\QQ$-cartier divisors on $Y$ are $f$-big.)
  
  Thus, the pair  $(Y, f_\ast^{-1} \Delta  + (1-\epsilon) \sum E_i)$ satisfies the hypotheses of 
   \cite[Theorem 1.2]{BCHM}, and we conclude \[
\bigoplus_{m\geq 0}
f_\ast \cO_Y( \lfloor m (
 K_Y+ f_\ast^{-1}\Delta +(1-\epsilon) \sum E_i ) \rfloor) 
 \]
 is a finitely generated $\cO_X$-algebra. Additionally, so is 
  \[
\bigoplus_{m\geq 0}
f_\ast \cO_Y( \lfloor m (
(a-1) E_Y+
  \sum_{i=1}^{r} (\ld{E_i}{X}{\Delta} -\epsilon ) E_i 
) \rfloor) 
 .\]
Since $\ld{E_i}{X}{\Delta} -\epsilon>0$ for all $i$, it follows 
\[
f_\ast \cO_Y( \lfloor m (
(a-1) E_W+
  \sum_{i=1}^{r} (\ld{E_i}{X}{\Delta} -\epsilon ) E_i ) \rfloor ) = f_\ast \cO_Y ( \lfloor m (a-1) E_W \rfloor ).\] After taking a proper Veronese of the previous graded ring, the result follows.
\end{proof}

\begin{lemma}\label{mldkk}
Let $E$ be a divisor over an affine klt variety $X$. The following are equivalent:
\begin{enumerate}
\item There exists an effective $\QQ$-divisor $\Delta$ such that $(X,\Delta)$ is klt  pair and
$\ld{E}{X}{\Delta}<1$. 
\item $k_E < \lct(\aeb)$. 
\end{enumerate}
\end{lemma} 
\begin{proof} We first show
(1) implies (2).
Assume the existence of such $\QQ$-divisor $\Delta$. Since $X$ is $\QQ$-factorial, we can choose $m\in \ZZ_{>0}$ such that $m\Delta$ is Cartier. Thus, $(X, \cO_X(-m \Delta)^{1/m})$ is a klt pair 
and 
$\ld{E}{X}{\cO_X(-m \Delta)^{1/m}} < 1$. By Proposition \ref{lctexp}, we conclude that $ k_E+1 - \lct(\aeb)< 1$ and the desired inequality follows. 
\\

We now show (2) implies (1). We first choose $m\in \ZZ_{>0}$ so that 
$k_E<m\cdot \lct(\fa_m^E)$, $\lct(\fa_m^E)<1$, and $\ord_E(\fa_m^E) = m$ (Lemma \ref{self}).  As described in \cite[Proposition 9.2.28]{MR2095472}, for a general element $f\in \fa_m^E$,
\[
\lct(f)= \lct(\fa_m^E)\hspace{.2 in} \text{ and } \hspace{.2 in} \ord_E(f)=m.\]

Now, set $\Delta= \lct(f) \{f =0 \}$, and note that $(X,\Delta)$ is log canonical by construction. 
Additionally, 
\[ 
\ld{E}{X}{\Delta} = k_E+1 -  \lct(f)  \ord_E(f) = k_E+1 -m \cdot \lct(\fa_m^E)<1 
,\]
where the last inequality follows from our assumption that $k_E<m \cdot \lct(\fa_m^E)$. While $(X,\Delta)$ is not klt, we note that $\Delta' = (1-\epsilon) \Delta$ satisfies the properties of (1) for $0<\epsilon \ll 1$. \end{proof}

\begin{proof}[Proof of Proposition \ref{proplctgeom}]  We first prove (i). The condition that $\aeb$ is finitely generated is local on $X$. Thus, it is sufficient to consider the case when $X$ is affine (which is necessary for the application of Lemma \ref{mldkk}). Now, it follows from Proposition \ref{mldk} and Lemma \ref{mldkk} that $\aeb$ is finitely generated.

We move on to (ii). By (i), $\aeb$ is finitely generated in degree. 
Therefore, (ii) follows from Proposition \ref{propblowup}. Furthermore, assume the $\aeb$ is finitely generated in degree $n$. As explained in the proof of Proposition \ref{propblowup}, $Y$ is the normalized blowup of $\fa_n^E$ and $\fa_n^E \cdot \cO_Y = \cO_Y(-nE_Y)$. 

We use the previous information to prove (iii). We first show
that $Y$ is $\QQ$-Gorenstein. Indeed, $K_Y= f^\ast K_X + K_{Y/X}$ is $\QQ$-Cartier. The divisor  $f^\ast K_X$ is $\QQ$-Cartier, since $K_X$ is $\QQ$-Cartier by assumption. The divisor $K_{Y/X}$ is $\QQ$-Cartier, since it is supported on $E_Y$ and $E_Y$ is $\QQ$-Cartier. 

Now, let $F$ be a divisor over $X$. We seek to compute $\ld{F}{Y}{\la E_Y}$ for $\la \geq 0$. 
We first compute $k_{F,Y}$ in terms of $k_{F,X}$.
Choose a projective birational morphism  $g:Z\to Y$ such that $Z$ is normal and  $F$ arises a prime divisor on $Z$. 
Since $K_{Z/X} = K_{Z/Y} + g^\ast K_{Y/X}$
and $K_{Y/X} = k_{E,X} E_Y$, we see
\[
k_{F,X} = k_{F,Y} + k_{E,X} \ord_F(E_Y)
.\]
We also have 
 \[\ord_F(E_Y)
 = \ord_F( \cO_{X_E} (-nE_Y) =
  \ord_F( \fa_n^E \cdot \cO_{Y} )  ,
 \]
 where the first equality is definitional and the second follows from the description of $Y$ as a blowup along $\fa_n^E$. 
We note
$\ord_F( \ab^E)= \ord_F(\fa_n^E)/n$, since $\aeb$ is finitely generated in degree $n$.

Finally, we compute
\[
\ld{F}{Y}{\la E_Y}=(
k_{F,X} - k_{E,X} \ord_F(E_Y)  +1) - \la \ord_F(E_Y) =
k_{F,X} - (k_{E,X} +\la) \ord_F( \aeb)
.\] 
Note that
\[
\lct(\ab^E) = \sup \{ \la \, \vert \, k_{F,X} +1 - \la \ord_F(\ab^E) > 0 \, \text{ for all } F \text{ over } X\}.\]
Since $\lct(\aeb)-k_E>0$, we see $Y$ is klt and  $\lct(Y,E_Y)= \lct(\aeb)-k_E$. 
\end{proof}

\begin{definition}\label{pltdef}
. Let $g:Z\to X$ be a projective birational morphism with exactly one irreducible exceptional divisor, say $S$. We say that $g:(Z,S) \to X$ is a \emph{plt blow-up} if $(Z,S)$ is plt and $-(K_Z+S)$ is $f$-ample \cite{plt}. When $g(S)$ is a closed point, $S$ is called a \emph{Koll{\'a}r component} \cite{Xu}.   

We note that such plt blowups were constructed in  \cite[Lemma 1]{Xu}.  Inspired by the work of Xu and others, we interpret these plt blowups in our framework. 
\end{definition}

\begin{proposition}\label{plt}
If  $E$ is a divisor over a klt variety $X$ with
\[
k_E+1 < \frac{ k_F+1}{ \ord_F(\aeb)},
\]
for all divisors $F$ over $X$ not equal to $E$,
then $(Y,E_Y)$ is plt where $Y:= \Proj( \oplus_{m\geq 0} \fa_m^E)$ and $E_Y$ is the prime divisor on $Y$ identified with $E$. Additionally, when $\codim(c_X(E))\geq 2$, $Y \to X$ is a plt blowup.

\end{proposition}

\begin{remark}
The condition that $k_E+1 < \frac{ k_F+1}{ \ord_F(\aeb)}$ for all $F \neq E$ implies that $\lct(\aeb) = k_E+1$. The converse does not hold.
\end{remark}

\begin{proof}
 Since $\lct(\aeb) = k_E+1$, Proposition \ref{proplctgeom} says that $(Y,E_Y)$ is log canonical. For a divisor $F\neq E$ over $Y$,
we have 
\[
\ld{F}{Y}{E_Y} =
k_{F,X} - (k_{E,X} +1) \ord_F( \aeb) \]
as in the proof of Proposition \ref{proplctgeom}. By our assumption, the latter value is $>0$. Thus, $(Y,E_Y)$ is plt.
 
We claim that  $-(K_Y+E_Y)$ is $f$-ample. Indeed,  $-(K_Y+E_Y)\sim_f -(K_{Y}-f^*(K_X) +E_Y) = -(k_E+1) E_Y$. Since $X$ is klt, $k_E+1>0$. Thus, $-(K_Y+E_Y)$ is $f$-linearly equivalent to a negative multiple of $E_Y$. Since $-E_Y$ is $f$-ample by Proposition \ref{propblowup}, we are done.
\end{proof}

\begin{lemma}
Let $X$ be a klt variety and $(X, \fa^\la)$ a pair
 with log resolution $f:Y\to X$. Let $\Delta_Y$  be the divisor on $Y$ such that $K_Y + \Delta_Y = f^\ast(K_X) + \la D$ where $\fa \cdot \cO_Y= \cO_Y(-D)$. 
 If $E$ is an exceptional divisor of $f$ and the coefficients of $\Delta_Y$ are  $\leq 1$ with equality precisely along $E$, then 
 \[
k_E+1 < \frac{ k_F+1}{ \ord_F(\aeb)}. 
\]
for all divisors $F$ over $X$ not equal to $E$.
\end{lemma}

\begin{proof}
Let $F$ be a divisor over $X$ with $F$ not equal to $E$. 
By \cite[Lemma 2.30]{KM}, we have $\ld{F}{X}{\fa^\la}=\ld{F}{Y}{\Delta_Y}$. 
By \cite[Lemma 2.30]{KM} and our hypotheses on the coefficients of $\Delta_Y$, we see $\ld{F}{X}{\fa^\la}=\ld{F}{Y}{\Delta_Y} >0$. 
By similar logic, $\ld{E}{X}{\fa^\la} =\ld{E}{Y}{\Delta_Y}=1$. Thus,
\begin{equation}\label{inequal}
k_F+1 > \la \ord_F(\fa) \hspace{0.5 in} \text{ and } \hspace{0.5 in}
k_E+1=\la \ord_E(\fa).
\end{equation}
Next, we note 
\[
\ord_F(\fa) \geq \ord_F(  \fa_{\ord_E(\fa)}^E) \geq  \ord_E( \fa ) \ord_E( \aeb)
,\]
where the first inequality comes from the inclusion $\fa \subseteq \fa_{\ord_E(\fa)}^E$ and the second from the definition of $ \ord_E( \aeb)$ as an infimum. Thus, the desired inequality follows.
\end{proof}

\begin{remark}
The previous lemma gives a condition for when the hypotheses of Proposition \ref{plt} hold. By the ``perturbation trick'' \cite[Proposition 3.2]{Ambrolocus}, if $\fa \subseteq \cO_X$ is a non-zero ideal on klt variety, then at least one of the divisors computing $\lct(\fa)$ satisfy the hypotheses of Proposition \ref{plt}. 

Thus,  if $\fa \subseteq \cO_X$ is a non-zero ideal on klt variety $X$, then at least of of the divisors computing $\lct(\fa)$ satisfy the conclusion of Proposition \ref{plt}. 
We note that a similar result was independently obtained by Kento Fujita in
\cite{Fujita}.

\end{remark}

\section{Connection with divisors that comput lct's and mld's}\label{sectionmld}

We proceed to prove  Theorem \ref{lct}. Recall that the theorem says that the value of $\lct(\aeb)$ determines whether $E$ computes a log canonical threshold.

\begin{proof}[Proof of Theorem  \ref{lct}]

The implication (1) implies (2) is trivial. Indeed, if $E$ compute $\lct(\fb)$, then $E$ computes $\lct(\fb_\bullet)$, where $\fb_\bullet$ is the trivial graded sequence defined by $\fb_m:= \fb^m$. 

Next, we show (2) implies (3). We follow an argument in \cite{JM}.
Assume $E$ computes $\lct(\ab)$, where $\ab$ is a graded sequence of ideals on $X$. 
By Lemma \ref{self}, we know that $\lct(\aeb) \leq k_E+1$. To prove that $\lct(\aeb)= k_E+1$, we must show that for any divisor  $F$ over $X$ we have $k_E+1 \leq (k_F+1)/ \ord_F(\aeb)$. 
Since $E$ computes $\lct(\ab)$, we have 
\[
\frac{k_E+1}{\ord_E(\ab) } = \lct(\ab) \leq \frac{k_F+1}{ \ord_F(\ab)}.\]
Thus, it is sufficient to show that $\ord_E(\ab) \cdot \ord_F(\aeb) \leq \ord_F(\ab)$.

We proceed to prove the previous inequality. 
By the definition of $\ord_E(\fa_m)$ as an infimum, we have $m\cdot \ord_E(\ab) \leq \ord_E(\fa_m)$ for all $m \in \NN$. Thus, 
$ \fa_m \subseteq \fa_{ \ord_E( \fa_m)}^E \subseteq \fa_{ \lfloor m\cdot \ord_E(\ab) \rfloor }^E$.
By the previous inclusion, $\ord_F( \fa_m) \geq \ord_F( \fa_{ \lfloor m\cdot \ord_E(\ab) \rfloor })$. After dividing by $m$ and taking infimums, the desired inequality follows.  \\

Lastly, we show (3) implies (1). Assume $\lct(\aeb)=k_E+1$. By Proposition \ref{proplctgeom}, the graded sequence $\aeb$ is finitely generated. Choose $n$ so that $\aeb$ is finitely generated in degree $n$. By our assumption,
\[
k_E+1=
\lct(\aeb)
=
n\cdot \lct(\fa_n^E)\]
Since $\ord_E(\fa^E_n)=n$, $E$ computes the log canonical threshold of $\fa_n^E$. \end{proof}

\begin{proposition}\label{jmp}
Let $E$ be a divisor over a klt  variety $X$. The divisor $E$ computes an asymptotic log canonical threshold if and only if $\lct(\aeb)$. 
\end{proposition}

\begin{proof}
Assume $E$ computes $\lct(\ab)$. As in (1) implies (2) of the previous proposition, 
we have that 
\end{proof}

Next, we seek to determine which divisors over a variety $X$ can compute minimal log discrepancies.

\begin{question}
Let $E$ be a divisor over $X$ such that $E$ computes a minimal log discrepancy.
\begin{enumerate}
\item Does it follow that $\aeb$ is finitely generated?
\item Does this imply that $\lct(\aeb)>k_E$ or $=k_E+1$?
\end{enumerate}
\end{question}
By Proposition \ref{lct} and Proposition \ref{proplctgeom}, the answer is yes to both if we replace the words ``minimal log discrepancy'' with ``log canonical threshold."
Also, note that, by Theorem \ref{proplctgeom}, an affirmative answer to the second question implies an affirmative
answer to the first. 

Let's consider a divisor $E$ such that $\lct(\aeb)<k_E+1$ and ask whether or not it can compute a minimal log discrepancy. 
We know that there must exist some divisor $F$ such that $\frac{k_F+1}{\ord_F(\aeb)}<k_E+1$. 
Additionally, it might be that $c_X(E)=c_X(F)$ and $k_F\leq k_E$. 
The proposition below shows that if such an $F$ exists, then $E$ cannot compute a minimal log discrepancy. 

\begin{proposition}\label{nomld} Let $E,F$ be divisors over a klt variety $X$ with
\begin{enumerate}
 \item $c_X(E)=c_X(F)$,
\item $k_F \leq k_E$, and
\item $ \displaystyle \frac{k_F+1}{\ord_F(\aeb)} <k_E+1.$
\end{enumerate}
Then for any log canonical pair $(X, \fb^\la)$ where $c_X(E) \subset \Cosupp(\fb)$ 
and $\la>0$
\[
\ld{F}{X}{\fb^\la} <\ld{E}{X}{\fb^\la}.
\]
\end{proposition}

Before proving the above proposition, we define the log discrepancy of a divisor $F$ over $X$ along 
$\ab^\la$, where $\ab$ is a graded sequence of ideals on $X$ with $S(\ab)\neq \emptyset$ and $\la\in \RR_{\geq 0}$, as 
\[
\ld{F}{X}{\ab^\la} \coloneqq  k_F+1-\lambda \ord_F(\ab)\]
and prove the following lemma. 

\begin{lemma}
Let $E,F$ be divisors over $X$ satisfying the conditions in Proposition \ref{nomld}. 
Then, 
\[
\ld{F}{X}{ \left( \aeb \right) ^\la} < \ld{E}{X}{ \left( \aeb \right)^\la }
\]
for all $\la\in \RR_{>0}$ such that $\ld{F}{X}{ \left( \aeb \right) ^\la }\geq 0$. 
\end{lemma}
\begin{proof}
Note that $\ld{ F}{X}{ \left( \aeb \right) ^\la }$ and $\ld{ F}{X}{ \left( \aeb \right) ^\la }$ are
real valued linear functions in $\la$.

When $\la=0$, we compare the values of the two functions:
\[
\ld{ E}{X}{\left( \aeb \right) ^0 }=k_E+1 \geq k_F+1 =\ld{ F}{X}{ \left( \aeb \right) ^0 }
.\]
Set $\la_F$ and $\la_E$ to be the values of $\la$ so that $\ld{F}{X}{ \left( \aeb \right) ^\la }=0$ and $\ld{ F}{X}{ \left( \aeb \right) ^\la }=0$, respectively. (The existence of $\la_F$ relies on $\ord_F(\aeb)>0$, which is implied by assumption (3).)  Note that 
\[
\la_F = \frac{k_E+1}{\ord_F(\aeb)} < k_E+1 =   \frac{k_E+1}{\ord_E(\aeb)}= \la_E
\]
where the inequality comes from our assumption. By analyzing these linear functions, we see that the
desired inequality holds. 
\end{proof}

\begin{proof}[Proof of Proposition \ref{nomld}]
Let $E,F$ and $(X,\fb^\la)$ satisfy the hypotheses of the proposition. We have 
\[
\ld{F}{X}{\fb^\la} \leq \ld{ F}{X}{ \left( \fa_{\ord_E(\fb)}^E \right) ^\la } \leq
\ld{ F}{X}{ \left( \aeb \right) ^{\la \ord_E(\fb)}}
<
\ld{ E}{X}{ \left( \aeb \right) ^{\la \ord_E(\fb)} }
=\ld{E}{X}{\fb^\la}\]
where the first inequality follows from $\fb \subset \fa_{\ord_E(\fb)}^E$, the second from the 
definition of the asymptotic order of vanishing as an infimum, the third from the previous lemma, 
and the last from Lemma \ref{self}. 
\end{proof}

\section{Divisors computing mlds on surfaces}\label{sectionsurf}

Before proving Theorem \ref{surf}, we prove the following lemma and proposition. Much of the content of
the following lemma can be found in \cite{Lip}.
\begin{lemma}\label{contract}
Let $x\in X$ be a point on a surface with at worst rational singularities and $E$ a divisor over $X$ with $c_X(E)= \{x\}$. 
If
\[
p : Y \to X
\]
is a projective birational morphism and $Y$ is a smooth variety that contains $E$ as a prime divisor,
then there is an $m\in \ZZ_{>0}$ such that $\aeb$ is finitely generated in degree $m$ and $\aem \cdot \cO_Y= \cO(-D)$, for some divisor $D$ on $Y$ with $\Supp(D) \subseteq \Exc(p)$.
\end{lemma}

\begin{proof}
The statement in the theorem is local in a neighborhood of $x\in X$. Thus, it is sufficient to consider the case when $X$ is a projective surface and $X$ is smooth or has an isolated singularity at $X$. Now, we consider the intersection form on the curves in $\Exc(f)$. By the Hodge Index Theorem, the intersection form on $\Exc(f)$ is negative definite. Thus, we may define
a $\QQ$-divisor $\check{E}$ with support on $\Exc(f)$ such that 
\[
\check{E}  \cdot C = 
\begin{cases}
    1, & \text{if $C=E$}.\\
    0, & \text{if $C\neq E$ and $C\subset \Exc(f)$}.
  \end{cases}
\]
Since $\check{E}$ intersects non-negatively with all exceptional curves
of $f$, $\check{E}$ is $f$-nef. By \cite[Theorem 12.1]{Lip}, $\check{E}$ is also $f$-base point free
(we are using that $X$ has rational singularities).  
Thus, $\check{E}$ gives rise to a fiber space $h$, over $X$, that contracts all curves in $\Exc(f)$
not equal to $E$. 
\
\[
\begin{tikzcd}
Y \arrow{rr}{h}\arrow{rd}{p} & &Z  \arrow{ld}{j} \\
  & X  &
  \end{tikzcd}
  \]

 The only exceptional divisor of the map $j$ (labeled in the above diagram) corresponds to $E$. 
 Since $j$ is a projective birational morphism with exactly one exceptional divisor corresponding to $E$, we conclude that $\aeb$ is finitely generated by Theorem \ref{tfg}. (To apply the previous theorem, we must have that $X$ is $\QQ$-factorial. By \cite[Proposition 17.1]{Lip}, rational surfaces singularities are $\QQ$-factorial.)
Additionally, $Z$ is isomorphic to the blowup of $X$ along $\fa_m^E$ for $m$ divisible enough (see the proof of Proposition \ref{tfg}).
For such an $m$, $\aem \cdot \cO_{Z}= \cO_{Z}(-mE)$ and $\aem\cdot \cO_Y = \cO_{Z}(-mE)\cdot \cO_Y = \cO_Y(h^\ast (mE))$. 
\end{proof}

\begin{proposition} \label{duv}
Let $X$ be a surface, $x\in X$ a smooth point or a du Val Singularity, and $(X,\fb^\la)$ a pair such that $x\in \Cosupp(\fb)$ and $\la>0$. If there exists a divisor $E$ that computes $\mld_x(X, \fb^\la)$, then $\lct(\aeb)=k_E+1$. 
\end{proposition}

\begin{proof}
We first consider the following maps
\[
Y \overset{f}{\longrightarrow} X' \overset{g}{\longrightarrow}X ,\]
where $g:X'\to X$ is the minimal resolution of $X$ and $f:Y \to X'$ be the map achieved by repeatedly blowing up the center of $E$ until the center is a prime divisor.
If $x\in X$ is a smooth point, then $X'=X$ and $g$ is the identity. Similarly, if $E$ corresponds to a prime divisor on $X'$, then $Y=X'$ and $f$ is the identity.

 Let $E_1,\ldots, E_s$ denote the exceptional divisors of $g\circ f$. 
Note that if $E_i$ is not contracted by $f$, then $f(E_i)$ is an exceptional divisor of $g$ and $k_{E_i}=0$. This follows from the definition of a du Val singularity.

By the previous lemma, there exists an $m\in \NN$ so that $\aeb$ is finitely generated in degree $m$ and $Y\to X$ is a log resolution of $\fa_m^E$. By Proposition \ref{tlct}, 
\[
\lct(\aeb) = \min_{i} \frac{k_{E_i} +1}{\ord_{E_i}(\aeb)}
.\]
To show that $\lct(\aeb)=k_E+1$, it is sufficient to show that
\[
k_E+1
\leq
\frac{k_{E_i}+1}{\ord_{E_i}(\aeb) }
\]
for all $1\leq i \leq s$. We claim that that $k_{E_i}  \leq k_E$ for all $1\leq i \leq s$.  By Proposition \ref{nomld}, if $k_E
>
\frac{k_{E_i}+1}{\ord_{E_i}(\aeb)}$, then $E$ cannot compute $\mld_x(X,\fb^\la)$. A contradiction. Thus, it suffices to prove our claim that
$k_{E_i}  \leq k_E$ for all $1\leq i \leq s$.

\noindent{\bf Case 1}: The map $f$ is the identity. 

In this case, $x\in X$ is a du Val singularity and $E$ corresponds to a prime divisor on the minimal resolution of $X$.
Thus, $k_{E_i} =0$ for all $1\leq i \leq s$ and $k_E=0$. We are done. 

\noindent {\bf Case 2}: The map $f$ is not the identity.

Let $E_r,E_{r+1},\ldots, E_s$ denote the exceptional divisors of $g\circ f$ that are contracted by $f$. Since these divisors arose via a sequence of blowups, we assume that the divisors are labelled in the order in which they arose. Thus, $E_s$ and $E$ are equivalent divisors over $X$.

If $1\leq i <r$, then $E_i$ is a an exceptional divisor of $g$ and, as stated before, $k_{E_i}=0$. 
If $r\leq j \leq s$, then $E_j$ either arose as the blowup of a point lying on a single exceptional divisor or the intersection of two such exceptional divisors. Thus, we have that either
\[
k_{E_j} = k_{E_{j-1}}+1 \text{ or } k_{E_j} = k_{E_{j-1}}+ k_{E_q} +1
\]  
for some $q<j-1$.
We see that $0\leq k_{E_{i-1}}\leq k_{E_i}$ for all $1<i\leq s$. Since $k_{E_s}=k_E$, we are done. \end{proof}

\begin{proof}[Proof of Theorem \ref{surf}]

Let $X$ be a smooth surface and $E$ a divisor over $X$ computing $\mld_x(X, \fb^\la)$ where $(X,\fb^\la)$ is a pair.

If $x\notin \Cosupp(\fb)$ or $\la=0$, then $\la \ord_F(\fb)=0$ for all divisors $F$ over $X$ with $\{x\}= c_X(F)$. Thus, 
\[
\mld_x(X, \fb^\la) =\min \{ k_E+1 \vert \, E \text{ is a divisor with } c_X(E)=\{x\} \} =2
\, \]
and the minimum is achieved by exactly one divisor over $X$, the divisor corresponding to the exceptional divisor of the blowup of $X$ at $\{x\}$.  Note that this divisor also computes a log canonical threshold, namely $\lct(\fm_x)$, where $\fm_x$ is the ideal of functions vanishing at $x$.  

If $x\in \Cosupp(\fb)$ and $\la>0$, then we apply the previous proposition to see that $\lct(\aeb)=k_E+1$. By Theorem \ref{lct}, $E$ must also compute a log canonical threshold. 
\end{proof}

\begin{remark}\label{duvr}
We can prove an analogus result for du Val Singularities. If $x\in X$ is a du Val singularity and $E$ is a divisor over $X$ computing $\mld_x(X,\fb^\la)$ where $x\in \Cosupp(\fb)$ and $\la>0$, then $E$ computes an log canonical threshold. This, does not include the case of $\mld_x(X,\cO_X)$. 

In fact, if $x\in X$ is an $E_7$ du Val singularity, then every exceptional divisor on the minimal resolution of $x\in X$ computes $\mld_x(X,\cO_X)$ (all exceptional divisors on the minimal resolution have log discrepancy $1$.) However, not of all of these exceptional divisors compute log canonical thresholds. 
\end{remark}

\section{Examples}\label{examples}

Below, we compute $\lct(\aeb)$ for a few examples of divisors over smooth varieties. 

\begin{example} For the simplest possible example consider $\Af^n$ when $n>1$ and let $E$ be the exceptional divisor  of the blowup of $\Af^n$ at the origin. Then, $\fa_m^E = \mathfrak{m}_0^m$, where $\mathfrak{m}_0 \subset \cO_{\Af^n}$ is the ideal of functions that vanish at the origin. Thus, $\lct(\aeb) = k_E+1 = n$ and is computed by $E$. 
\end{example}

\begin{example}
Let $X$ be a smooth surface and   $ Y \to X$ the composition of $r\geq 3$ point blowups resulting in the following dual diagram of exceptional curves
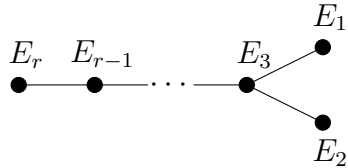
\begin{center}
\begin{tikzpicture}

	\draw (1,0) -- (2.7,0);
	\draw (3.3, 0) -- (4,0);
	\draw (4,0) -- (5,-.5);
	\draw (4,0) -- (5,.5);
	
	\draw[fill=black] (1,0) circle(.1);
	\draw[fill=black] (2,0) circle(.1);
	\draw[fill=black] (4,0) circle(.1);
	\draw[fill=black] (5,-.5) circle(.1);
	\draw[fill=black] (5,.5) circle(.1);
	
	\node at (3,0) {$\cdots$};
	\node at (5.1,.9) {$E_1$};
	\node at (5.1,-.9) {$E_2$};
	\node at (4.1,.4) {$E_3$};
	\node at (2.1,.4) {$E_{r-1}$};
	 \node at (1.1,.4) {$E_{r}$};
\end{tikzpicture}
\end{center}
where $E_i$ denotes the strict transform of the exceptional divisor arising from the $i$-th blowup. To understand $\fa_\bullet^{E_r}$ for $r\geq 3$, we compute the intersection form on $\Exc(Y)$ and apply the argument in the proof of Lemma \ref{contract}.  We find that $\fa_\bullet^{E_r}$ is finitely generated in degree $r+3$ 
and 
\[
\fa_{r+3}^{E_r} \cdot \cO_Y  = \cO_Y\left(-2E_1 -3E_2 - \sum_{i=3}^r (i+3) E_i \right)
.\]
Since $k_{E_1}=1$, $k_{E_2} = 2$, and $k_{E_i}= i+1$ for $3 \leq i \leq r$, we see that
\[
\lct({\fa_\bullet}^{E_r} ) = \frac{k_{E_3}+1}{ \ord_{E_3} ({\fa_\bullet}^{E_r})} = \frac{r+2}{6/(r+3) } = 6\left( \frac{r+2}{r+3} \right). \]

Thus, we see two behaviors.  We have that $\lct( \fa^{E_r}_\bullet) =k_{E_r}+1$ when $r=3$, but this is not the case when $r>3$.  The divisor $E_3$ computes a log canonical threshold, while $E_r$ does not for $r>3$ (Proposition \ref{lct}).

Additionally, the divisor $E_r$ does not compute a minimal log discrepancy for $r>3$ (Theorem \ref{surf}). We may view $E_3$ as preventing $E_r$ from computing a minimal log discrepancy. For any log canonical pair $(X,\fb^\la)$ and $r>3$, we have that 
\[
\ld{E_r}{X}{\fb^\la} > \ld{E_3}{X}{\fb^\la}
\]
by Proposition \ref{nomld}.

\end{example}

\begin{example} \label{nfgex} For an example of a divisor over a smooth variety with a non-finitely generated graded sequence of ideals, we look to \cite{Kur}. 
A divisor $E$ over  $\Af^4$ is constructed by the composition of maps
\[
E \subset Y \overset{f}{\longrightarrow} B_0\Af^4 \overset{g}{\longrightarrow} \Af^4
.\]
The map $g$ is the blowup of $\Af^4$ at the origin with exceptional divisor $F$. The map $f$ is the blowup of a specially chosen smooth curve inside $F \simeq \PP^3$ and $E$ is the exceptional divisor of this map. 
Since $\aeb$ is not finitely generated, there does not exist a morphism $Y\to X_E$ over $X$ that contracts $F$ (Theorem \ref{tfg}). Intuitively, $F$  ``obstructs'' the finite generation of $\aeb$. Additionally, 
 \[
 \lct(\aeb) \leq \frac{k_F+1}{\ord_F(\aeb)} \leq 5 < 6 = k_E+1
  \,.\]
  where the first inequality comes from the definition of $\lct(\aeb)$ and the second inequality from the fact that $\ord_F(\aeb) \geq  4/5$ and $k_F=3$.
The divisor $E$ cannot compute a log canonical threshold (Proposition \ref{lct}) or minimal log discrepancy (Proposition \ref{nomld}).

\end{example}

\end{document}